 \RequirePackage{}
\documentclass[smallextended]{svjour3}       
\smartqed  
\usepackage [latin1]{inputenc}

\usepackage{graphicx}
\usepackage{epstopdf}
\usepackage{algorithm} 
\usepackage{algorithmic} 
\usepackage[colorlinks,bookmarksopen,bookmarksnumbered,citecolor=blue, linkcolor=blue, urlcolor=blue]{hyperref}

\usepackage{amsmath}
\usepackage{enumerate}
\usepackage{multirow}
\usepackage{booktabs}
\usepackage{amssymb}
 \usepackage{mathrsfs}
\usepackage{graphicx}
\usepackage{subfigure}
\usepackage{rotfloat}
\usepackage{rotating} 
\usepackage{epstopdf,ntheorem}
\usepackage{amsfonts}
\usepackage{caption}
\newtheorem{remark}{Remark}
\newtheorem{theorem}{Theorem}
\newtheorem{lemma}{Lemma}
\begin{document}
	
	\title{
		The convergence rate of the accelerated proximal gradient
		algorithm for Multiobjective Optimization is  faster than 
		$O(1/k^2)$ 
	}
	
	
	\author{Jie Zhang        \and  Xinmin Yang\textsuperscript{*}
	}
	
	
	\institute{School of Science, Chongqing University of Posts and Telecommunications, Chonngqing, China (Jie Zhang)
		\\
		\email{zjieabc@163.com} (Jie Zhang)            \\
		School of Mathematics Science, Chongqing Normal University, Chongqing, China(Xinmin Yang) \\
 	\email{xmyang@cqnu.edu.cn}(Xinmin Yang)\\
		*  Corresponding author
	}
	
	\date{Received: date / Accepted: date}

	\maketitle
	\begin{abstract}	 
		In this paper, we propose a fast  proximal gradient  algorithm for multiobjective optimization,   it is proved that the convergence rate  of the  accelerated algorithm for multiobjective optimization  developed by  Tanabe et al. can be improved from  $O(1/k^2)$ to $o(1/k^2)$ by introducing  different extrapolation term  $\frac{k-1}{k+\alpha-1}$ with $\alpha>3$.  Further, we establish the inexact version of the proposed algorithm when the error term is additive, which owns the same convergence rate.
		At last, the efficiency of the proposed algorithm is verified on some numerical experiments.
	\end{abstract}
\section{Introduction}
 We mainly focus on the multiobjective optimization,  
\begin{align}\label{0.02}
	\min\limits_{x\in\mathbb{R}^n} F(x)
\end{align}
with $F:\mathbb{R}^n\to(\mathbb{R}\cup\{\infty\})^m$  and
$F:=(F_1,\dots,F_m)^T$ taking the form
\begin{align}\label{0.04}
 F_i(x):=f_i(x)+g_i(x),\, i=1,2,\dots,m,
\end{align}
where $f_i:\mathbb{R}^n\to\mathbb{R}$ is a convex and continuously differentiable function and $g_i:\mathbb{R}^n\to(-\infty,+\infty] $ is a closed, proper and convex function. Suppose the gradient $\nabla f_i(x)$ of $ f_i(x) $ in (\ref{0.04}) is Lipschitz continuous with $L_i$ and let
\[ L:= \max_{i=1,\dots,m} L_i,
\]

Multiobjective optimization refers to the process of simultaneously minimizing (or maximizing) multiple objective functions while taking into account any relevant constraints. In most cases, it is not possible to find a single point that minimizes all objective functions at the same time, making the concept of Pareto optimality indispensable. A point is deemed Pareto optimal or efficient if there is no other point with equal or lower objective function values, and at least one objective function value is strictly lower. The applications of multi-objective optimization span across diverse domains, and for an extensive list of examples, please refer to \cite{Beck,Boyd}.

Recently, the  descent methods for multiobjective optimization problems, 

Various first-order methods, including steepest descent types,\cite{Fliege2}, projected gradient methods,\cite{Fliege2}, proximal point methods,\cite{Fukuda}, conjugate gradient methods,\cite{Lucambio}, Barzilai-Borwein descent methods,\cite{ChenJian}, and conditional gradient methods,\cite{chenw}, have garnered significant attention and extensive research efforts. These methods are widely studied for optimization tasks. 
In addition to first-order approaches, second-order methods such as the Newton method,\cite{Fliege1}, trust region methods,\cite{Carrizo}, and Newton-type proximal gradient methods,\cite{Ansary} have been explored. These second-order methods leverage information about the curvature of the objective function, making them more efficient for certain optimization problems but potentially more computationally expensive.

When $m=1$, the multiobjective optimization problem (MOP) (\ref{0.02}) reduces to the single objective case,
\begin{align}\label{0.03}
	\min\limits_{x\in\mathbb{R}^n} F(x):=f(x)+g(x).
\end{align}
The proximal gradient (PG) algorithm is a popular method to solve the composition optimization problem especially when the proximal operator of $g$ has an explicit solution. However, the convergence rate $O(1/k)$ of the proximal gradient method is slow. Lots of researchers devote to improve the efficiency of the PG algorithm. In 2009, Beck and Teboulle \cite{Beck}
developed the remarkble Fast Iterative Shrinkage Thresholding Algorithm (FISTA) for the convex case. The convergence rate $O(1/k^2)$ of the objective function is obtained.  Moreover, Chambolle and Dossal \cite{Chambolle} gave the modification version of the FISTA and further  obtained the sequence convergence of proposed algorithm.It is worth noting that Attouch and Peypouquet \cite{Attouch} introduced an extrapolation technique, specifically $\frac{k-1}{k+\alpha-1}$ with $\alpha>3$, into the proximal gradient method for solving (\ref{0.03}). This innovation led to an enhanced convergence rate, improving from $O(1/k^2)$ to $o(1/k^2)$.

The accelerated algorithm mentioned above
does not leverage additional second-order information at a higher computational cost. Despite this, it achieves an enhanced convergence rate.
Inspired by the benefit of acceleration algorithms on single-objective scenarios,  the researchers are increasingly focusing on the study of acceleration algorithms in multi-objective optimization problems. Recently, Tanabe et al. \cite{Tanabe1} extends the remarkable FISTA to the multiobjective case and the convergence rate $O(1/k^2)$ characterized by merit function \cite{Tanabe} is obtained which improved the proximal gradient method for MOP \cite{Tanabe3}.  
 Nishimura et al. \cite{Nishimura}  established the monotonicity
version of the multiobjective FISTA. Moreover, Tanabe et al.  \cite{Tanabe4} generalizes the multiobjective FISTA by introducing some 
hyperparameters which is also a generalization even in the single-objective case and preserves the same convergence rate $O(1/k^2)$
 as  multiobjective FISTA.  Besides, it is proved that the   iterative sequences is convergent. 
  
Based on the fast algorithms in single-objective case     proposed by Attouch \cite{Attouch}, we introduce the extrapolation parameter $\frac{k-1}{k+\alpha-1}$ with $\alpha>3$ into the multiobjective proximal gradient algorithm. The convergence rate  described by the merit function achives $o(1/k^2)$ which improves the 
multiobjective FISTA \cite{Tanabe1}. Along with we get the sequence convergence of the iterations. Additionaly,  the inexact version of the proposed algorithm is also got when the errors satisfying some additivity and the convergence rate remains unchanged.
 
The organization of the paper is as follows. Some notations and definitions are given in Sect. 2. Sect. 3 presents the algorithm framework of fast proximal gradient algorithm for multi-objective optimization. 
 Sect. 4  focus on proving the convergence of the proposed algorithm. Numerical results
are presented in Sect. 5, which demonstrates that the proposed algorithm  exhibits rapid convergence effects. Towards the conclusion of the paper, several key conclusions are presented.

\section{Preliminary Results}
For convenience, we present some notations and definitions which used in the overall paper.
 Denote $\mathbb{R}^n$ be the $n$-dimensional space, and $\mathbb{R}^n_{+}:=
 \{v\in\mathbb{R}^n|v_i\ge 0,i=1,2,\dots,n\}$; $\triangle^n$ be the standard
 simplex in $\mathbb{R}^n$ with 
 \[  \triangle^n:=\left\{\lambda\in\mathbb{R}^n_{+}\bigg|\lambda_i\ge 0, \sum_{i=1}^{n}\lambda_i=1.   \right\}
 \] 
The Euclidean norm defined in $\mathbb{R}^n$ is  
$\Vert u\Vert:=(\sum_{i=1}^{n} u^2_i)^{1/2}$ for $u\in\mathbb{R}^n$.

 The partial orders induced by $\mathbb{R}^n_{+}$ implies that 
for any $v^1,v^2\in\mathbb{R}^n_{+}$, if $v^2-v^1\in\mathbb{R}^n_{+} $, then 
$ v^1\leq v^2 $  and if $v^2-v^1\in {\rm int}\mathbb{R}^n_{+} $, then $ v^1< v^2 $.

For a closed, proper and convex function
 $h:\mathbb{R}^n\to\mathbb{R}\cup\{\infty\}$, the Moreau envelope of $h$ defined by
 \[  \mathcal{M}_h(x):=\min_{y\in\mathbb{R}^n} \left\{ h(y)+\frac{1}{2}\Vert x-y\Vert^2 \right\}
 \]
 The unique solution of the above problem is called the 
 proximal operator of $h$, which is denoted as
 \[{\rm prox}_h(x):
 =\arg\min_{y\in\mathbb{R}^n}\left\{ h(y)+\frac{1}{2}\Vert x-y \Vert^2 \right\}.
 \]
 \begin{lemma}{\rm \cite{Beck}}\label{lemma2}
  If $h$ is a proper closed and convex function, the Moreau 
  envelope $ \mathcal{M}_h$ is 1-Lipschitz continuous and its gradient takes the following form,
  \[ \nabla\mathcal{M}_h(x):=x- {\rm prox}_h(x).  
  \] 
 \end{lemma}

Next we recall the optimality concept for the multiobjective 
optimization problem (\ref{0.02}): $x^*\in\mathbb{R}^n$ is weakly pareto optimal if there no exist $x\in\mathbb{R}^n$ such that  $F(x)<F(x^*)$. Denote the set of weakly pareto optimal solutions be $X^*$.
 
The merit function $u_0: \mathbb{R}^n\to\mathbb{R}\cup\{\infty \}$ given in \cite{Tanabe} takes the following form
\begin{align}\label{0.05}
 u_0(x):= \sup\limits_{z\in\mathbb{R}^n}\min_{i=1,\dots,m}[F_i(x)-F_i(z)],
 \end{align}
 and proves that $u_0$ is a merit function in the Pareto sense.
\begin{lemma} {\rm\cite{Tanabe}}\label{lemma3}
	Let $u_0$ be given as (\ref{0.05}), then $u_0(x)\ge 0$, $x\in\mathbb{R}^n$, and $x$ is the weakly 
	Pareto optimal for (\ref{0.02}) if and only if $u_0(x)=0$.
\end{lemma}
	\section{The fast proximal gradient algorithm for multi-objective optimization}
In this section, we propose the fast proximal gradient algorithm for multi-objective optimization. The subproblem  has the same form  in \cite{Tanabe}, in detail, for $x\in dom F,\, \ell>L$, 
\begin{align}\label{0.1}
\min_{z\in\mathbb{R}^n} \phi_{\ell}(z;x,y) 
\end{align}
where 
\begin{align}\label{0.3}
 \phi_{\ell}(z;x,y):
=\max\limits_{i=1,\dots,m}\left\{\langle\nabla f_i(y),z-y\rangle +g_i(z)+f_i(y)-F_i(x) \right\}+\frac{\ell}{2}\Vert z-y\Vert^2.
\end{align}
Define
\begin{align}\label{0.4}
p_{\ell}(x,y):=
\arg\min\limits_{z\in\mathbb{R}^n}\phi_{\ell}(z;x,y),\quad  \theta_{\ell}(x,y):=\min_{z\in\mathbb{R}^n}\phi_{\ell}(z;x,y)
\end{align}

The optimality condition of (\ref{0.1}) means that
 there exists 
$\eta(x,y)\in\partial g(p_{\ell}(x,y))$ and the Lagrange multiplier $\lambda(x,y)\in\mathbb{R}^m$ satisfying
\begin{align}
	\sum_{i=1}^{m} &\lambda_i(x,y)\{\nabla f_i(y)+\eta_i(x,y)  \}=-\ell(p_{\ell}(x,y)-y) \label{0.5}	\\	
		&\lambda(x,y)\in\triangle^m,  \lambda_i(x,y)=0, i\notin \mathcal{I}(x,y),\label{0.6}
\end{align} 
where $\triangle^m$ represents the standard simplex and 
\[
\mathcal{I}(x,y):=\arg\max\limits_{i=1,\cdots,m}\left\{\langle\nabla f_i(y),z-y\rangle +g_i(z)+f_i(y)-F_i(x) \right\}+\frac{\ell}{2}\Vert z-y\Vert^2.
\]

\begin{lemma}\rm\cite{Tanabe1}
	Suppose $p_{\ell}(x,y)$ is defined as in (\ref{0.4}),
	then $y\in\mathbb{R}^n$ is the weakly Pareto optimal for (\ref{0.02}) if and only if $p_{\ell}(x,y)=y$ for some $x\in\mathbb{R}^n$.
\end{lemma}

The result implies that it is allowable to use  $\Vert p^{acc}_{\ell}(x,y)-y   \Vert_{\infty}$ as the stopping criterion of the corresponding fast proximal gradient algorithm for MOP. 
The specific algorithm framework is as follows. 
	\begin{algorithm}[hptb]
		\caption{The Fast Proximal Gradient Algorithm for Multi-objective Optimization } 
		\label{alg1} 
		\begin{algorithmic}[1] 
			\STATE 
					
		\textbf{Initialization:} Give a starting point $x^0=y^1\in dom F$, $\alpha>3$, $\varepsilon>0$.			
		  Set $k=0$.
		  \WHILE{ $\Vert p_{\ell}(x^{k-1},y^k)-y^k\Vert\ge \varepsilon$}
			\STATE 
		\noindent 
	\[ x^k\to p_{\ell}(x^{k-1},y^k) ,\]
			\STATE   
			\begin{align}\label{0.2}
			y^{k+1}=x^k+\frac{k-1}{k+\alpha-1}(x^k-x^{k-1}),
			\end{align}  			  
				\STATE 	\noindent $k\to k+1$.

	 	\ENDWHILE
		\end{algorithmic}
	\end{algorithm}

\section{The convergence rate analysis of FPGMOP for (\ref{0.02}).}
 Given a fixed weakly pareto solution $x^*\in\mathbb{R}^n$,   define the auxiliary sequence
  \begin{align}\label{0.01}
   \mathscr{E}_{k+1}:=\frac{2(k+\alpha-1)^2}{\ell(\alpha-1)}\sigma_{k+1}+(\alpha-1)\Vert u^{k+1}-x^* \Vert^2.
   \end{align}	
where $\sigma_{k+1}:=\sigma_{k+1}(x^*)=\min\limits_{i=1,\cdots,m}\left\{F_i(x^k)-F_i(x^*)\right\}$, $u^{k+1}:=\frac{k+\alpha-1}{\alpha-1}x^{k+1}-\frac{k}{\alpha-1}x^{k}$.


 \begin{lemma}\label{lemma1}{\rm\cite{Tanabe1}}
 Suppse $\{x^k\}$ and $\{y^k\}$ are the sequences genetated by Algorithm\ref{alg1}, it holds that for any $k\ge0$

 \begin{align}\label{1.1}
 \sigma_{k+1}(z)\leq \frac{-\ell}{2}\left\{2\langle x^{k+1}-y^{k+1},y^{k+1}-z \rangle +\Vert x^{k+1}-y^{k+1}\Vert^2  \right\}-\frac{\ell-L}{2}\Vert x^{k+1}-y^{k+1} \Vert^2 ,
 \end{align}
 and 
\begin{align}\label{1.2}
 \sigma_{k}(z)-\sigma_{k+1}(z)\ge
 \frac{\ell}{2}\left\{2\langle x^{k+1}-y^{k+1},y^{k+1}-x^k \rangle +\Vert x^{k+1}-y^{k+1}\Vert^2  \right\}+\frac{\ell-L}{2}\Vert x^{k+1}-y^{k+1} \Vert^2,
 \end{align} 
where $\sigma_k(z)=\min\limits_{i=1,\cdots,m}\left\{F_i(x^k)-F_i(z)\right\}$.
  \end{lemma}

\begin{theorem}\label{theorem1}
	Suppose $\{x^k\}$ and $\{y^k\}$ be the sequences genetated by Algorithm\ref{alg1}, for any $z\in\mathbb{R}^n$, $k\ge0$ and  $\alpha>3$,
	it holds that
	\begin{enumerate}[\rm(i)]
		\item  $\mathscr{E}_{k+1}\leq M:=\frac{2(\alpha-1)}{l}\Vert x^0-x^* \Vert^2.$	 
		\item  $\sum_{k=0}^{+\infty}(k+\alpha-1)\sigma_{k}(z) \leq\frac{\ell M(\alpha-1) }{\alpha-3};$
		\item  $\sum_{k=0}^{+\infty}k\Vert x^{k+1}-x^k\Vert^2\leq \frac{3\ell(\alpha-1)M}{(\alpha-3)};$
		\item   
		$\lim_{{k\to+\infty}} k^2\Vert x^{k+1}-x^k\Vert^2+(k+1)^2\sigma_k \quad \text{\rm exists};
		$
	 
	\item  $ u_0(x^k)=o(1/k^2) \quad \text{\rm and} \quad
		 \Vert x^{k+1}-x^k \Vert=o(1/k) .$
		
		\end{enumerate} 
	
\end{theorem} 
{\it{Proof.}}
(i)
	Accoring to Lemma \ref{lemma1},  multiplying (\ref{1.1}) by $\frac{\alpha-1}{k+\alpha-1}$ and
	 (\ref{1.2}) by $\frac{k}{k+\alpha-1} $ and adding them together, it holds that
	 \begin{align*} \sigma_{k+1}(z)-\frac{k}{k+\alpha-1}\sigma_{k}(z)\leq
	 \ell\left\{\langle y^{k+1}-x^{k+1},y^{k+1}
	 -\frac{k}{k+\alpha-1}x^k-\frac{\alpha-1}{k+\alpha-1}z  \rangle \right\}	  \\ 
	  -\frac{\ell}{2}\Vert x^{k+1}-y^{k+1}\Vert^2  -   \frac{\ell-L}{2}\Vert x^{k+1}-y^{k+1} \Vert^2,
	 \end{align*}
	 It is noted that 
	 \begin{align*}
	 &-\Vert x^{k+1}-y^{k+1}\Vert^2 
	 +2\langle y^{k+1}-x^{k+1},  y^{k+1}
	 -\frac{k}{k+\alpha-1}x^k-\frac{\alpha-1}{k+\alpha-1}z   \rangle \\
	 =&-\Vert x^{k+1}-\frac{k}{k+\alpha-1}x^k -\frac{\alpha-1}{k+\alpha-1}z \Vert^2
	 +\Vert y^{k+1}-\frac{k}{k+\alpha-1}x^k -\frac{\alpha-1}{k+\alpha-1}z \Vert^2 \\
	 =&-(\frac{\alpha-1}{k+\alpha-1})^2(\Vert u^{k+1}-z \Vert^2-\Vert u^k-z\Vert^2 )	 
	 	 \end{align*}
 	 where $u^{k}:=\frac{k+\alpha-2}{\alpha-1}x^k-\frac{k-1}{\alpha-1}x^{k-1}$.  	 
 	 Then we get
 	 \[ 
 	 \sigma_{k+1}(z)-\frac{k}{k+\alpha-1}\sigma_{k}(z)
 	 \leq -\frac{\ell}{2}(\frac{\alpha-1}{k+\alpha-1})^2(\Vert u^{k+1}-z \Vert^2-\Vert u^k-z\Vert^2 )-\frac{\ell-L}{2}\Vert x^{k+1}-y^{k+1} \Vert^2,
 	 \]
 	 
Multiplying the above inequality by $\frac{2(k+\alpha-1)^2}{\ell(\alpha-1)}$, it yields that 
\[
 \mathscr{E}_{k+1}\leq \frac{2k(k+\alpha-1)}{\ell(\alpha-1)}\sigma_{k}(z)
+(\alpha-1)\Vert u^k-z \Vert^2-\frac{(\ell-L)(k+\alpha-1)^2}{\ell(\alpha-1)}\Vert x^{k+1}-y^{k+1} \Vert^2.
\]

 Since $k(k+\alpha-1)<(k+\alpha-2)^2-(\alpha-3)(k+\alpha-1)$, we get
 \begin{align}\label{1.3}
 \mathscr{E}_{k+1}+\frac{2(\alpha-3)(k+\alpha-1)}{\ell(\alpha-1)}\sigma_{k}(z)
 \leq \mathscr{E}_{k}  -\frac{(\ell-L)(k+\alpha-1)^2}{\ell(\alpha-1)}\Vert x^{k+1}-y^{k+1} \Vert^2\leq \mathscr{E}_{k},
 \end{align} 
 Then the sequence $\{\mathscr{E}_{k} \}$ is nonincreasing.
 Moreover, from (\ref{0.01}) and $\sigma_k(z)\ge0$, it yields that $\mathscr{E}_{k}$ is lower-bounded  
  and $\lim_{k\to+\infty} \mathscr{E}_{k} $ exists and 
   ${\mathscr{E}_{k} }\leq {\mathscr{E}_{1}} $ for $k\ge0$ .
  
  From (\ref{1.1}), we have 
  \begin{align*} 
 	\sigma_{1}&\leq \frac{-\ell}{2}\left\{2\langle x^{1}-y^{1},y^{1}-x^* \rangle +\Vert x^{1}-y^{1}\Vert^2  \right\}-\frac{\ell-L}{2}\Vert x^{1}-y^{1} \Vert^2 ,\\
 	&\leq \frac{-\ell}{2}\left\{2\langle x^{1}-y^{1},y^{1}-x^* \rangle +\Vert x^{1}-y^{1}\Vert^2  \right\}\\
 	&\leq \frac{\ell}{2}(-\Vert x^1-x^*\Vert^2+\Vert y^1-x^*\Vert^2)=\frac{\ell}{2}(-\Vert x^1-x^*\Vert^2+\Vert x^0-x^*\Vert^2)
 \end{align*}
Combining this with  $\mathscr{E}_1=\frac{2(\alpha-1)}{l}\sigma_1
+(\alpha-1)\Vert x^1-x^* \Vert^2$, we get  
\[ {\mathscr{E}_{1} }\leq \frac{2(\alpha-1)}{l} \Vert x^0-x^* \Vert^2,    \]
i.e.,for any $k\ge 0$, it holds that
\[  \mathscr{E}_{k+1}\leq \frac{2(\alpha-1)}{l}\Vert x^0-x^* \Vert^2.  \]
So the result (i) hold.
 
(ii)
 Summing up the inequality (\ref{1.3}) from $k=1$ to $k=K$, we get 
\begin{align}\label{1.4}
 \mathscr{E}_{K+1}+\sum_{k=1}^{K} \frac{(\alpha-3)(k+\alpha-1)}{\ell(\alpha-1)}\sigma_{k}(z)\leq \mathscr{E}_{1}.
 \end{align}

From (\ref{1.4}), we have 

\[ \sum_{k=0}^{K} \frac{(\alpha-3)(k+\alpha-1)}{\ell(\alpha-1)}\sigma_{k}(z)\leq \mathscr{E}_{1}- \mathscr{E}_{K+1}.
\] 
Letting $K\to+\infty$, due to $ \mathscr{E}_{K+1}\leq\mathscr{E}_{1}$ and $\lim_{k\to\infty}\mathscr{E}_{k+1}$ exists, we get
\[ \sum_{k=0}^{+\infty}(k+\alpha-1)\sigma_{k}(z)\leq \frac{\ell\mathscr{E}_{1}(\alpha-1) }{\alpha-3},\]
so the result (ii) holds.

 (iii) According to (\ref{1.2}) and the Pythagoras relation 
 \begin{align}\label{1.55}
   \Vert b-a \Vert^2+2\langle b-a,a-c \rangle=\Vert b-c \Vert^2-\Vert a-c \Vert^2 
 \end{align}
  with $(a,b,c)=(y^{k+1},x^{k+1},x^k)$,
  we have  
 \begin{align}\label{1.5}
 	\sigma_{k+1}(z)-\sigma_k(z)&\leq  \frac{\ell}{2}\left\{-\Vert x^{k+1}-x^k\Vert^2+\Vert y^{k+1}-x^k\Vert^2 \right\} 
 	-\frac{\ell-L}{2}\Vert x^{k+1}-y^k\Vert^2  \notag\\
 	&=\frac{\ell}{2}\left\{-\Vert x^{k+1}-x^k\Vert^2
 	+(\frac{k-1}{k+\alpha-1})^2
 	\Vert x^{k}-x^{k-1}\Vert^2 \right\} \notag\\
 	&\quad-\frac{\ell-L}{2}\Vert x^{k+1}-x^k-\frac{k-1}{k+\alpha-1}(x^k-x^{k-1})\Vert^2  \notag\\
 	&= -\frac{2\ell-L}{2}\Vert x^{k+1}-x^k \Vert^2 +\frac{L}{2}(\frac{k-1}{k+\alpha-1})^2\Vert x^{k}-x^{k-1}\Vert^2  \notag\\
 	 &\quad 
 	 + (\ell-L)\langle x^{k+1}-x^k,\frac{k-1}{k+\alpha-1}(x^k-x^{k-1}) \rangle
 	  \notag\\
 	 &\leq-\frac{\ell}{2}\Vert x^{k+1}-x^k\Vert^2+\frac{\ell}{2}(\frac{k-1}{k+\alpha-1})^2\Vert x^{k}-x^{k-1}\Vert^2,  
 \end{align}
where the last inequaltiy holds due to the Caucy-Schwarz inequality. 

 The subsequent analysis is similar to  (17) in \cite{Attouch}, 
  for completeness, a detailed derivation is given below. 
 From (\ref{1.5}) and $k+\alpha-1\ge k+1$,  it holds that
 \[ (k+1)^2d_k-(k-1)^2d_{k-1}\leq (k+1)^2(	\sigma_{k}(z)-\sigma_{k+1}(z)),
 \]
   where $d_k=\frac{\ell}{2}\Vert x^{k+1}-x^k \Vert^2$.
    Since
   \begin{align*}
   (k+1)^2( \sigma_{k}(z)-\sigma_{k+1}(z) )
   &=k^2\sigma_{k}(z)-(k+1)^2\sigma_{k+1}(z)+(2k+1) \sigma_{k}(z) \notag\\  
   &\leq k^2\sigma_{k}(z)-(k+1)^2\sigma_{k+1}(z) +3k\sigma_{k}(z),
   \end{align*} 
 we get  
   \begin{align*}
   	2kd_k+k^2d_k-(k-1)^2d_{k-1}\leq &(k+1)^2d_k-(k-1)^2d_{k-1}\notag\\
   	\leq &(k+1)^2(	\sigma_{k}(z)-\sigma_{k+1}(z))\notag\\
   	\leq &k^2	\sigma_{k}(z)-(k+1)^2\sigma_{k+1}(z)+3k\sigma_{k}(z).
   	   	 \end{align*}
      	 Summing the inequality from $k=1$ to $K$, we have
      	  \[  K^2d_K+2\sum_{k=0}^{K}kd_{k}\leq \sigma_1+\frac{3\ell(\alpha-1)\mathscr{E}_1}{2(\alpha-3)},  \quad K\in\mathbb{N},
 \]
  and so
 \[ \sum_{k=0}^{+\infty} k\Vert x^{k+1}-x^{k}\Vert^2\leq\frac{3\ell(\alpha-1)\mathscr{E}_0}{(\alpha-3)}.
 \]
 
 (iv) Similar to the analysis of Lemma 2 in \cite{Attouch}, the desired results can be obtained.
 
  (v) From (\ref{0.3}) and the result (iv), it can be easily got. 
  
 In the following, we are ready to prove the convergence of the iterations. Now we first recall two  lemmas which will be used in the later. 
 \begin{lemma}{\rm\cite{Opial}}\label{lemma4}
 	(Opial's lemma) Let $S\subset\mathbb{R}^n$ be a nonempty set, if the sequence $\{z^k\}$ satisfies the following
 		\begin{enumerate}[\rm(i)]
 			\item ${\rm lim}_{k\to+\infty}\Vert z_k-z\Vert$ exists for every $z\in S$,
 			\item every limit point of the sequence $\{z_k\}$ is in the set $S$,  	
 	 	\end{enumerate}
  		\noindent then the sequence $\{z^k \}$ convergence to a point in $S$. 
 \end{lemma}
  \begin{lemma}\label{lemma5}{\rm\cite{Attouch}}
Suppose $\alpha>3$ and  nonegative sequences $\{a_k\}$ and $\{w_k\}$ satisfy
\[  a_{k+1}\leq \frac{k-1}{k+\alpha-1}a_k+w_k,
\]
and $\sum_{k=1}^{\infty}kw_k<\infty,$ then $\sum_{k=1}^{\infty} a_k<\infty. $
\end{lemma}
\begin{theorem}
	Suppose $\{x^k\}$ is the sequence generated by   algorithm \ref{alg1},  then the sequence gerenrated by the algorithm \ref{alg1} is weakly Pareto optimal for (\ref{0.02}). 
\end{theorem}
{\it{Proof.}} Based on the Lemma \ref{lemma3} and Theorem \ref{theorem1}, we just need to prove that the sequence $\{x^k \}$ generated by the algorithm \ref{alg1} convergences.
Since $\sigma_{k+1}(x^*)\ge 0$, $x^*\in X^*$,
from (\ref{1.1}), we get 
\begin{align*}
 2\langle y^{k+1}-x^{k+1}, y^{k+1}-x^* \rangle
 -\Vert x^{k+1}-y^{k+1}\Vert^2\ge 0,
\end{align*}
that is 
\[ \Vert y^{k+1}-x^* \Vert^2-\Vert x^{k+1}-x^* \Vert^2\ge 0.
\]
Then 
\begin{align*}
&	\Vert x^{k+1}-x^* \Vert^2\notag\\
	 \leq & \Vert y^{k+1}-x^* \Vert^2 \notag\\
	= &\Vert x^{k}-x^* \Vert^2 + (\frac{k-1}{k+\alpha-1})^2   \Vert x^{k}-x^{k-1} \Vert^2 + 2(\frac{k-1}{k+\alpha-1})\langle  x^{k}-x^*,x^{k}-x^{k-1}   \rangle  \notag\\
		=&\Vert x^{k}-x^* \Vert^2+ ((\frac{k-1}{k+\alpha-1})^2+\frac{k-1}{k+\alpha-1})\Vert x^{k}-x^{k-1} \Vert^2  \notag\\
		&\quad+ \frac{k-1}{k+\alpha-1}(\Vert x^{k}-x^* \Vert^2-\Vert x^{k-1}-x^* \Vert^2) 	\notag\\
		\leq& \Vert x^{k}-x^* \Vert^2 +2\Vert x^{k}-x^{k-1} \Vert^2++ \frac{k-1}{k+\alpha-1}(\Vert x^{k}-x^* \Vert^2-\Vert x^{k-1}-x^* \Vert^2)
\end{align*}
Let $h_k:=\Vert x^{k}-x^*\Vert^2$, then we have
\[ (h_{k+1}-h_k)_+\leq (\frac{k-1}{k+\alpha-1})(h_k-h_{k-1})_{+}+ 2\Vert x^{k}-x^{k-1} \Vert^2,
\]
According to Theorem\ref{theorem1}(iii) and Lemma \ref{lemma5}, we get
$\sum_{k=1}^{\infty} (h_{k+1}-h_k)_+<\infty,$ and further we  have
$\lim_{k\to\infty} h_k$ exists by the nonegative of $\{h_k\}$.
\begin{remark}
Recently, 
it is observed that the fast multiobjective gradient methods with Nesterov acceleration also proposed through the inertial gradient-like dynamical systems \cite{Sonntag1,Sonntag2}. 
In \cite{Sonntag1},  the authors considered the problem (\ref{0.02}) with $h(x)=0$ and established the inertial first order
method for MOP by an discretization of the differential equation. Further, they proved the convergence rate characterised by the merit function is  $O(1/k^2)$ just when the extrapolation parameter 
$\beta_k=\frac{k-1}{k+2}$ which is the special case of (\ref{0.2});   
 this algorithm  has the similar convergence behaviors to the accelerated algorithm  given by \cite {Tanabe1} which has been discussed in \cite{Sonntag1}. However, the convergence rate of the accelerated proximal gradient algorithms when $\alpha>3$ has not been obtained both in the  discrete case \cite{Sonntag1} and  the continuous case \cite{Sonntag2}. This question also was raised as an open question in the Section 9 of \cite{Sonntag1}. Fortunately, based on the  Lemma 6.3, Lemma 6.5 in \cite{Sonntag1} and motivated by the Theorem \ref{theorem1}, we also can obtain that  the convergence rate $o(1/k^2)$ also can be obtained when $\alpha>3$ in the fast proximal gradient algorithm for MOP \cite{Sonntag1} which solved their open question 
 \cite{Sonntag1}.

\end{remark}

 \section{Stability} 
 In this part, we consider the inexact version of the Algorithm \ref{alg1}, the inexact error  mainly comes from to the 
 calculation of the $\nabla f_i(x)$, in detail, the subproblem in (\ref{0.1}) takes the following form, 
 \begin{align}\label{2.1}
 	{\hat\phi}_{\ell}(z;x,y,\epsilon):
 	=\max\limits_{i=1,\dots,m}\left\{\langle\nabla f_i(y)+\epsilon,z-y\rangle +g_i(z)+f_i(y)-F_i(x) \right\}+\frac{L}{2}\Vert z-y\Vert^2,
 \end{align}
 where $\epsilon\in\mathbb{R}^n$ is the unknown error. 
 The corresponding optimality condition yields that there exists 
 ${\hat\eta}(x,y,\epsilon)\in\partial g(\hat p_{\ell}(x,y,\epsilon))$ and the Lagrange multiplier ${\hat\lambda}(x,y,\epsilon)\in\mathbb{R}^m$ satisfying
 \begin{align}
 	\sum_{i=1}^{m}&{\hat \lambda}_i(x,y,\epsilon)\{\nabla f_i(y)+\hat\eta_i(x,y,\epsilon)+\epsilon \}=-L(\hat p_{\ell}(x,y,\epsilon)-y) ,\label{2.2}\\ 	
 &\hat\lambda(x,y,\epsilon)\in\triangle^m,  \hat\lambda_j(x,y,\epsilon)=0, j\notin \mathcal{I}(x,y,\epsilon),\label{2.3}
 \end{align}
 where $\triangle^m$ represents the standard simplex and 
 \[
 \mathcal{I}(x,y,\epsilon):=\arg\max\limits_{i=1,\cdots,m}\left\{\langle\nabla f_i(y)+\epsilon,z-y\rangle +g_i(z)+f_i(y)-F_i(x) \right\}+\frac{\ell}{2}\Vert z-y\Vert^2.
 \]
 Define
 \begin{align}\label{2.04}
 \hat p_{\ell}(x,y,\epsilon):=
 	\arg\min\limits_{z\in\mathbb{R}^n}\hat\phi_{\ell}(z;x,y,\epsilon),\quad  \hat\theta_{\ell}(x,y,\epsilon):=\min_{z\in\mathbb{R}^n}\hat\phi_{\ell}(z;x,y,\epsilon)
 \end{align}
The algorithm framework is as follows:
 	\begin{algorithm}[H]
 	\caption{ The inexact Fast Proximal Gradient Algorithm for Multi-objective Optimization } 
 	\label{alg2} 
 	\begin{algorithmic}[1] 
 		\STATE 
 		
 		\textbf{Initialization:} Give a starting point $x^0=y^1\in\mathbb{R}^n$, $\alpha>3$.			
 		\WHILE{ $\Vert \hat p_{\ell}(x^{k-1},y^k,\epsilon^k)-y^k\Vert\ge \varepsilon$}
 		\STATE 
 		\noindent 
 		\[ x^k\to\hat p_{\ell}(x^{k-1},y^k,\epsilon^k) ,\]
 		\STATE   
 		\begin{align}\label{0.002}
 			y^{k+1}=x^k+\frac{k-1}{k+\alpha-1}(x^k-x^{k-1}),
 		\end{align}  			  
 		\STATE 	\noindent $k\to k+1$.

 		\ENDWHILE
 	\end{algorithmic}
 \end{algorithm}
 \begin{lemma} \label{lemma2.2}
 	Suppse $\{x^k\}$ and $\{y^k\}$ are the sequences genetated by Algorithm \ref{alg2}, it holds that
 \begin{align}\label{2.22}
  \sigma_{k+1}(z)\leq-\frac{L}{2}\left\{ \Vert x^{k+1}-z \Vert^2-\Vert y^{k+1}-z \Vert^2     \right\} 
  -\langle \epsilon^{k+1}, x^{k+1}-z  \rangle,
  \end{align}
  \begin{align}\label{2.222}
   \sigma_{k+1}(z)-\sigma_{k}(z)\leq-\frac{L}{2}\left\{ \Vert x^{k+1}-x^k \Vert^2-\Vert y^{k+1}-x^k \Vert^2     \right\} 
  -\langle \epsilon^{k+1}, x^{k+1}-x^k \rangle.
  \end{align} 
 {\it{Proof.}}
 By the analysis from \cite{Tanabe1} in Lemma 5.1, we get
\[
\sigma_{k+1}(z)\leq \sum_{i=1}^{m}\lambda_i(x^k,y^{k+1},\epsilon_{k+1})
\left\{\nabla f_i(y^{k+1})+\hat\eta_i(x^k,y^{k+1},\epsilon_{k+1}),x^{k+1}-z \right\}+
\frac{L}{2}\Vert x^{k+1}-y^{k+1}\Vert^2,
\]
 and then from (\ref{2.2}) and (\ref{2.3}), we have 
\begin{align*}
 \sigma_{k+1}(z)&\leq  \langle -L(x^{k+1}-y^{k+1})-\epsilon^{k+1}, x^{k+1}-z \rangle+\frac{L}{2} \Vert x^{k+1}-y^{k+1}  \Vert^2 \notag\\ 
 &=-\frac{L}{2}\left\{2\langle x^{k+1}-y^{k+1}+L^{-1}\epsilon^{k+1}, x^{k+1}-z \rangle -\Vert x^{k+1}-y^{k+1} \Vert^2  \right\} \notag\\ 
 &=-\frac{L}{2}\left\{2\langle x^{k+1}-y^{k+1}, y^{k+1}-z \rangle +\Vert x^{k+1}-y^{k+1} \Vert^2  \right\}-\langle \epsilon^{k+1}, x^{k+1}-z \rangle \notag\\ 
  &=  -\frac{L}{2} \left\{\Vert x^{k+1}-z\Vert^2-\Vert y^{k+1}-z\Vert^2  \right\}-\langle \epsilon^{k+1}, x^{k+1}-z \rangle, 
\end{align*}
so  we get the result (\ref{2.22}) and then we  infer the inequality  (\ref{2.22})  similarily as Lemma \ref{lemma1}.
 
 \end{lemma}
Before giving the main conclusions, we first present the following lemma,
\begin{lemma}{\rm\cite{Attouch2}}\label{lemma2.3}
	Suppose $\{a_k\}$ is a nonnegative sequence and satisfies 
	\[ a^2_k\leq c^2+\sum_{j=1}^k\beta_ja_j, \, k\in\mathbb{N},
	\]
	where $c>0$, $\beta_j>0$ and $\sum_{j=1}^{\infty}\beta_j<+\infty$, it has that
	$a_k\leq c+\sum_{j=1}^{\infty}\beta_j$.	  
\end{lemma}
\begin{theorem}
	Suppose $\{x^k\}$ is the sequence generated by the inexact version of the  algorithm \ref{alg2} and the errors satisfying $\sum_{k=0}^{\infty}(k+\alpha-1)\Vert \epsilon_k\Vert<+\infty$, then
	\begin{enumerate}[\rm(i)]
		\item $\sum_{k=0}^{\infty} (k+\alpha-1)\sigma_{k}<+\infty,$
		$\sum_{k=0}^{\infty}k\Vert x^{k+1}-x^{k} \Vert^2<\infty$;
		\item $\lim_{k\to+\infty}[(k-1)^2\Vert x^k-x^{k-1}\Vert^2+
		\frac{2}{L}(k+\alpha-2)^2\sigma_k]$ exists.
		\item the results in Theorem \ref{theorem1} {\rm(v)} hold.
		\end{enumerate} 
	\end{theorem}
 {\it{Proof.}}
  Define 
  \begin{align}\label{2.44}
  \hat{\mathscr{E}}_k={\mathscr{E}}_k +2L^{-1}\sum_{j=k}^{\infty}
  (j+\alpha-1)\langle \varepsilon^{k+1}_j, x^*-u^{k+1}_j  \rangle .
  \end{align}
  Similar to the analysis of Theorem \ref{theorem1}, we get 
   \begin{align*} 
  	\mathscr{E}_{k+1}+\frac{2(\alpha-3)(k+\alpha-1)}{L(\alpha-1)}\sigma_{k}(z)
  	\leq \mathscr{E}_{k}+2L^{-1}(k+\alpha-1)\langle\varepsilon^{k+1},x^*-u^{k+1} \rangle,
  \end{align*}  
Then we have 
\begin{align}\label{2.4} 
\hat{\mathscr{E}}_{k+1}+\frac{2(\alpha-3)(k+\alpha-1)}{L(\alpha-1)}\sigma_{k}(z)
\leq \hat{\mathscr{E}}_{k},
\end{align}
i.e., $\{\hat{\mathscr{E}}_{k} \}$ is non-increasing.
Note that $\sigma_{k}\ge 0 $, then
\begin{align*}
\Vert u^{k+1}-x^* \Vert^2\leq\frac{1}{\alpha-1}{\mathscr{E}}_{1}
+\frac{2L^{-1}}{\alpha-1}\sum_{j=1}^{\infty}(j+\alpha-1)\Vert\varepsilon^{k+1}_{j} \Vert\Vert u^{k+1}_j-x^* \Vert.
\end{align*}
By virtue of $\sum_{k=1}^{\infty}(k+\alpha-1)\Vert \varepsilon_k\Vert<+\infty$ and Lemma \ref{lemma2.3}, it yields that
\[
\Vert u^{k+1}-x^* \Vert\leq C, \,\text{with}\, C=\sqrt{\frac{1}{\alpha-1}{\mathscr{E}}_{1}}+\frac{2L^{-1}}{\alpha-1}\sum_{j=1}^{\infty}(j+\alpha-1)\Vert\varepsilon^{k+1}_{j} \Vert.
\]
Since $\mathscr{E}_{k}\ge 0 $,  we get
\[ \hat{\mathscr{E}_{k}}\ge -2CL^{-1}\sum_{j=k}^{\infty}
(j+\alpha-1)\Vert\varepsilon^{k+1}_j\Vert >-\infty.
\]
From the nonincreasing of $\{  \hat{\mathscr{E}}_{k}  \}$,  then we get $\lim_{k\to\infty} \hat{\mathscr{E}_{k}} $ exists.
Further, from (\ref{2.4}), we have
\begin{align}\label{2.444}
\sum_{k=1}^{\infty}(k+\alpha-1)\sigma_k(z)	<+\infty.
\end{align}

By (\ref{2.222}), we have
 \begin{align*} 
	\sigma_{k+1}(z)-\sigma_{k}(z)\leq-\frac{L}{2} \Vert x^{k+1}-x^k \Vert^2+\frac{L}{2}(\frac{k-1}{k+\alpha-1})^2\Vert x^{k}-x^{k-1} \Vert^2     
	-\langle \varepsilon^{k+1}, x^{k+1}-x^k \rangle.
\end{align*} 

Multiplying the above inequality by $(k+\alpha-1)^2$ on both sides, we get
\begin{align*}
	2(k+\alpha-1)^2(\sigma_{k+1}(z)-\sigma_{k}(z))
	&\leq -L(k+\alpha-1)^2\Vert x^{k+1}-x^k \Vert^2+L(k-1)^2\Vert x^{k}-x^{k-1} \Vert^2 \notag\\
	 &\quad	+(k+\alpha-1)^2\langle \varepsilon^{k+1}, x^{k}-x^{k+1} \rangle
\end{align*}
From
\[ (k+\alpha-1)^2\sigma_{k+1}(z)-(k+\alpha-2)^2\sigma_{k}(z)
=(k+\alpha-1)^2(\sigma_{k+1}(z)-\sigma_{k}(z))+(2k+2\alpha-3)\sigma_{k}(z),
\]
we obatin
\begin{align}\label{2.5}
		&2(k+\alpha-1)^2\sigma_{k+1}(z)-2(k+\alpha-2)^2\sigma_{k}(z)+
	L(k+\alpha-1)^2\Vert x^{k+1}-x^k \Vert^2 \notag\\
	&-L(k-1)^2\Vert x^{k}-x^{k-1} \Vert^2    
  \leq (k+\alpha-1)^2\langle \varepsilon^{k+1}, x^{k}-x^{k+1} \rangle +(2k+2\alpha-3)\sigma_{k}(z).
\end{align}

Summing the inequality (\ref{2.5}) from $k=1$ to $k=K$, by $\sigma_k(z)\ge 0$, we deduce that
\begin{align}\label{2.6}
	&	L(K+\alpha-1)^2\Vert x^{K+1}-x^K \Vert^2\notag\\
\leq&	L(K+\alpha-1)^2\Vert x^{K+1}-x^K \Vert^2 +L\sum_{k=1}^{K-1}((k+\alpha-1)^2-k^2)\Vert x^{k+1}-x^k \Vert^2\notag\\
	\leq& M^2+
	\sum_{k=1}^{K}(k+\alpha-1)^2\Vert\varepsilon^{k+1}\Vert\Vert  x^{k}-x^{k+1} \Vert ,
\end{align}
where $M:=\sqrt{2(\alpha-1)^2\sigma_1(z)+\sum_{k=1}^{\infty}(2k+2\alpha-3)\sigma_{k}(z)} $ and $M<\infty$ from (\ref{2.444}).

Further, from Lemma \ref{lemma2.3} with 
$a^k=\sqrt{L}(k+\alpha-1)\Vert x^k-x^{k+1} \Vert$ 
and $\beta_j=\sqrt{L}^{-1}(j+\alpha-1)\Vert\varepsilon_j\Vert $, we have
\begin{align}\label{2.7}
(k+\alpha-1)\Vert x^{k+1}-x^k\Vert\leq W:=C+\sum_{j=1}^{\infty}\beta_j<+\infty.
\end{align}
Combining the above inequality with (\ref{2.6}), by $\alpha>3 $, we get
\[
 4L\sum_{k=1}^{\infty}k\Vert x^{k+1}-x^k \Vert^2\leq M^2+2W\sum_{k=1}^{\infty}(k+\alpha-1)\Vert\varepsilon_k\Vert<\infty.
\]
so the conclusion (i) is verified.

(ii) Define 
\[\Lambda_k:= \lim_{k\to+\infty}[{L}(k-1)^2\Vert x^k-x^{k-1}\Vert^2+
 {2}(k+\alpha-2)^2\sigma_k]    \]
 From (\ref{2.5}) and (\ref{2.6}), we infer that 
 \[ \Lambda_{k+1}-\Lambda_{k}\leq (2k+2\alpha-3)\sigma_k+2W(k+\alpha-1)
 \Vert\varepsilon_k \Vert 
 \]
 Together this with (\ref{2.444}), we get  $\sum_{k=1}^{\infty}(k+\alpha-1)\Vert \varepsilon_k\Vert<+\infty$. Since $\Lambda_k\ge0$,  taking the positive part of the left inequality, then $\lim_{k\to\infty}\Lambda_k$ exists.
and we easily have that the results in Theorem \ref{theorem1} holds.

\section{Numerical experiments} 
This section provides a performance comparison of the proposed algorithm with two other algorithms for the multiobjective optimization: the proximal gradient method \cite{Tanabe3} (denoted as   SPG ) and the accelerated proximal gradient method \cite{Tanabe1} (denoted as   accelerated SPG) with the following  extrapolation parameter 
 \[\gamma_k=\frac{t_k-1}{t_k},\,  t_{k+1}=\sqrt{t^2_k+\frac{1}{4}}+\frac{1}{2},\, k\in\mathbb{N}.\] 

The original test problems encompass both constrained and unconstrained problems which are presented in Table \ref{table1}. In the case of constrained problems,   the indicator function $g_i$ is applied to represent the constraint. 
For unconstrained problems, there are two types being considered: $g_i=0$ and  $g_i = \Vert x - i + 1\Vert /(in)$, $i=1,2,\dots,m$.

The codes are running on a laptop with a processor model i7-11390H (3.4GHz) and 16GB of memory, with a 64-bit operating systema. Solving the original problem (\ref{0.3}) by solving the dual form of the subproblem based on \cite{Fliege1}. The dual problem is a problem with unit simplex constraints and it can be efficiently solved by the 
condition gradient method \cite{Beck}. The algorithm is stopped if   the subproblem is satisfied $\Vert  x^{k+1}-x^k\Vert\leq 10^{-11}$
or the maximum number of iterations reaches  $k_{\rm max}=2000$. 
For problem $FDS$, the parameter $\ell$ is determined by backtrackting procedutre with the initial value is 1, and the backtracking factor is 2. For the other problems, the parameter $\ell$ can be calculated accurately. The 100
  initial points are selected by uniformly and randomly within the given boundaries in \cite{Mita}.
  
   The average numerical results of the algorithms are displayed 
   in Table \ref{tab2}. It can be observed that the proposed algorithm 
   meet the stopping criteria within the less time and fewer iterations 
   for the test problems other than the FDS promblem with $g_i=l_1$   which  need more iterations but use the least time. 
   
\begin{figure}[htpb]
	\centering
	\begin{minipage}[t] {0.3\textwidth}
		\includegraphics[scale=0.03]{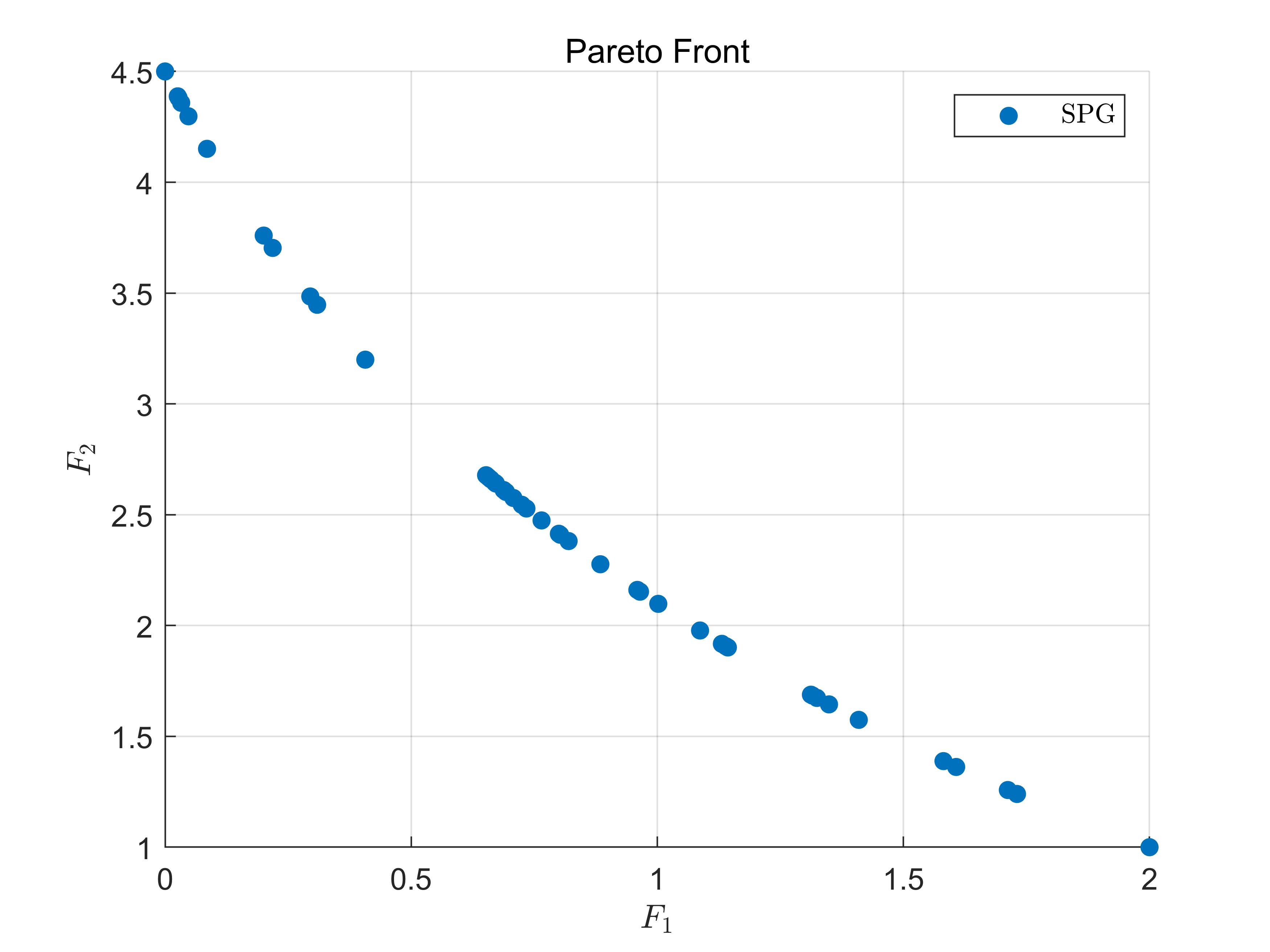}
	\end{minipage}	
	\centering
	\begin{minipage}[t] {0.3\textwidth}
		\includegraphics[scale=0.03]{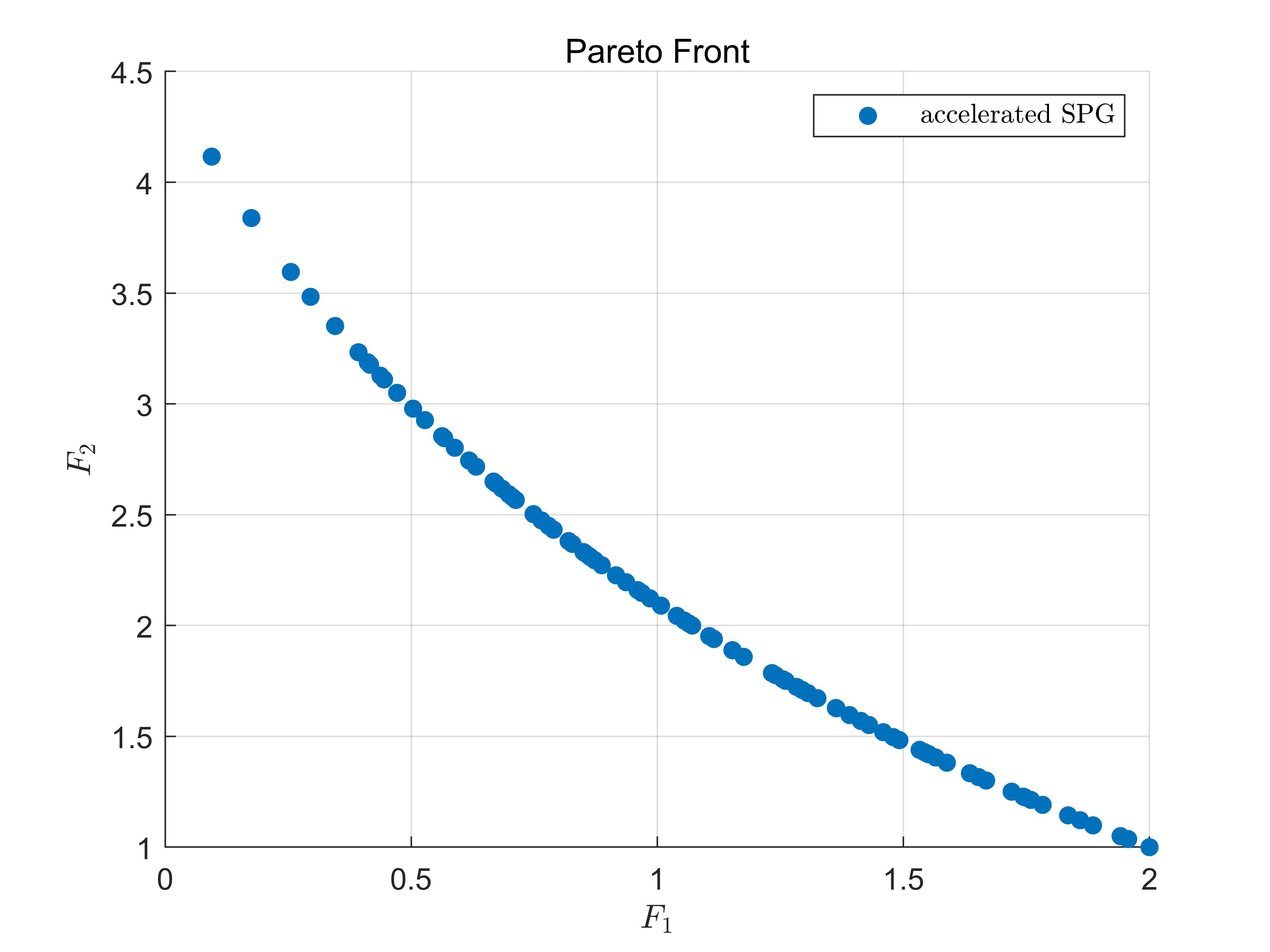}
	\end{minipage}
	\begin{minipage}[t] {0.3\textwidth}
		\includegraphics[scale=0.03]{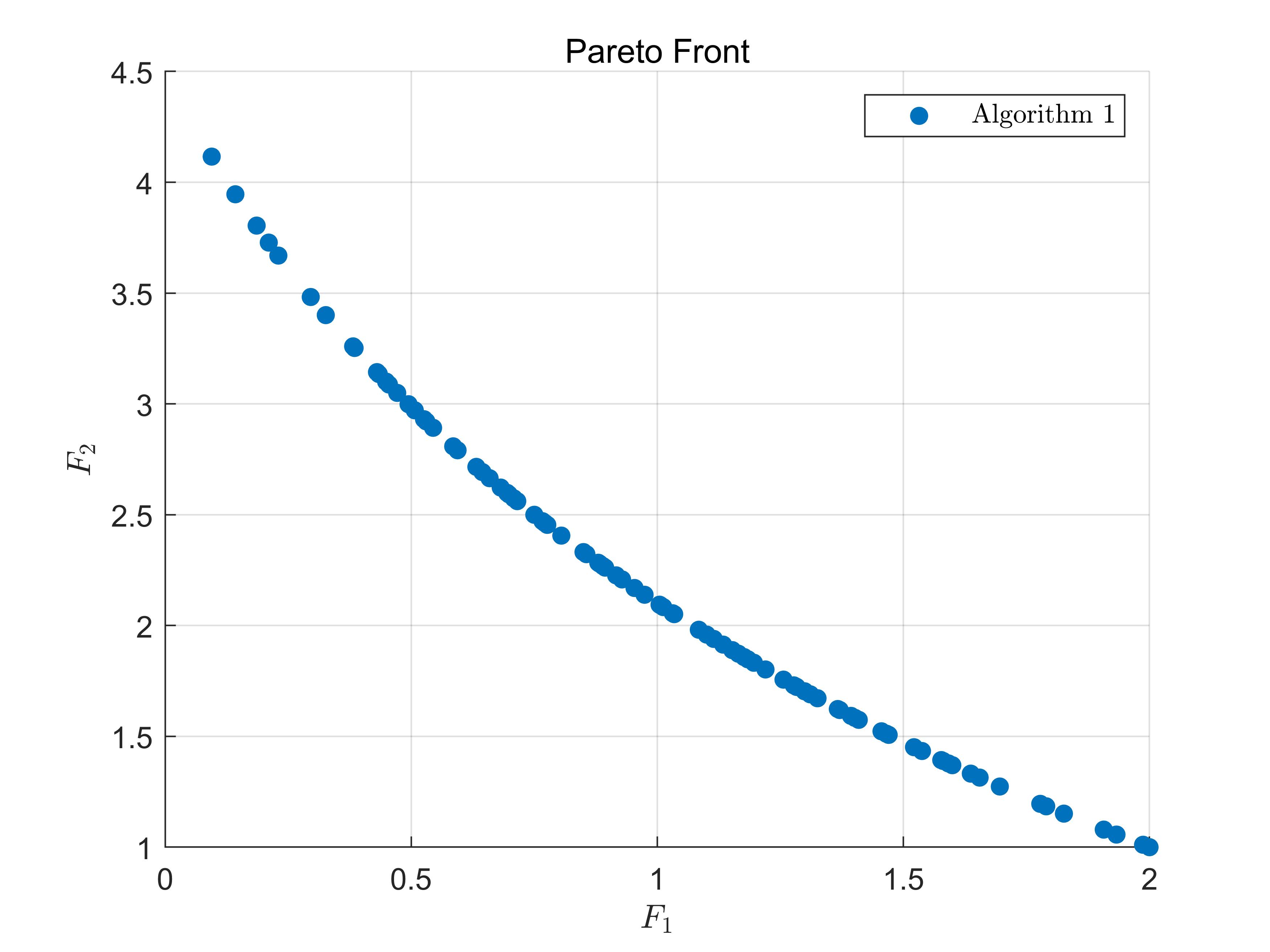}	
	\end{minipage}
	\centering\caption{The final Pareto front for JOS1 with $n=50$.}
	\label{fig1}
\end{figure}

\begin{figure}[htpb]
	\centering
	\begin{minipage}[t] {0.3\textwidth}
			\includegraphics[scale=0.03]{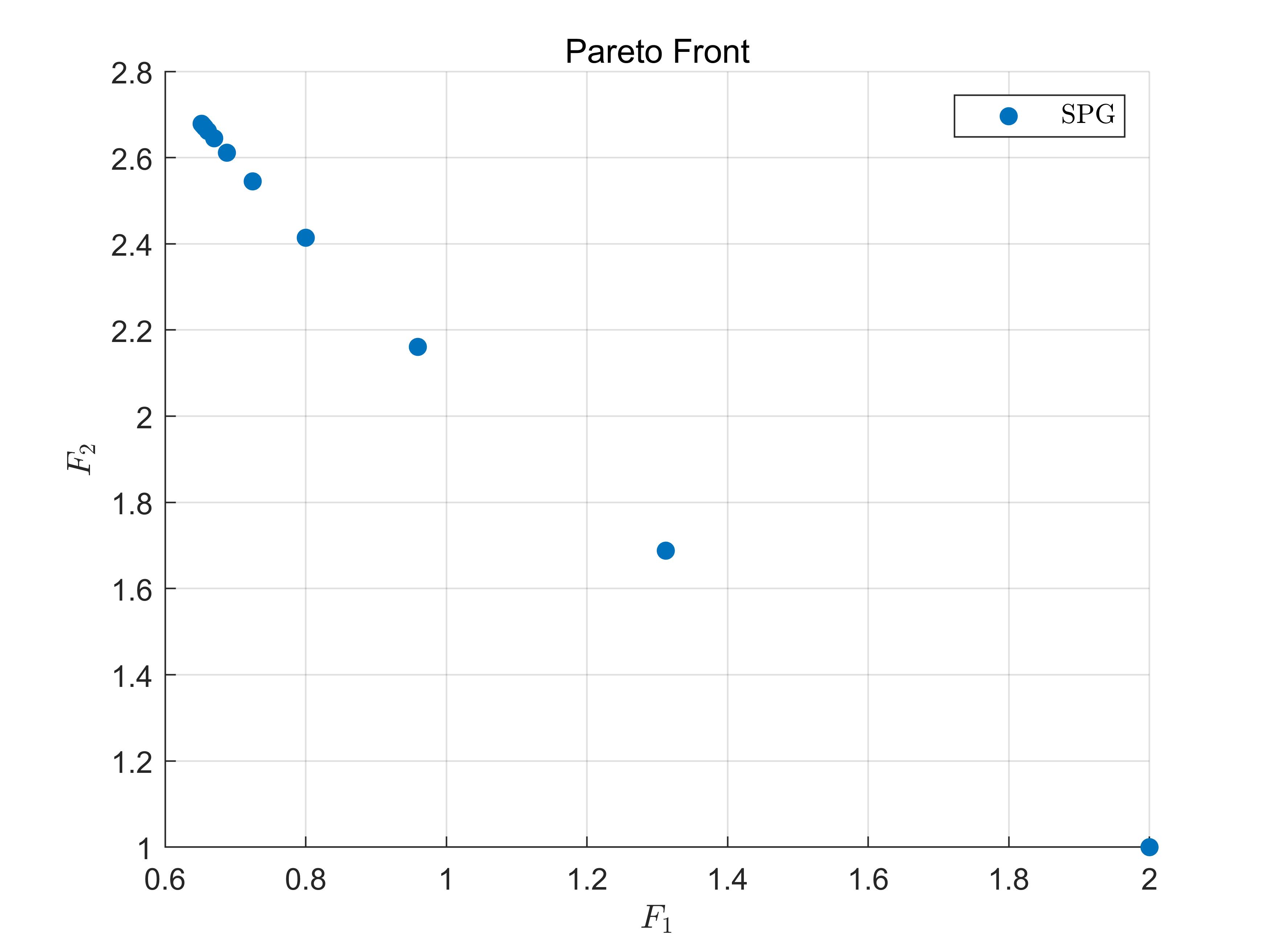}
	\end{minipage}	
	\centering
	\begin{minipage}[t] {0.3\textwidth}
			\includegraphics[scale=0.03]{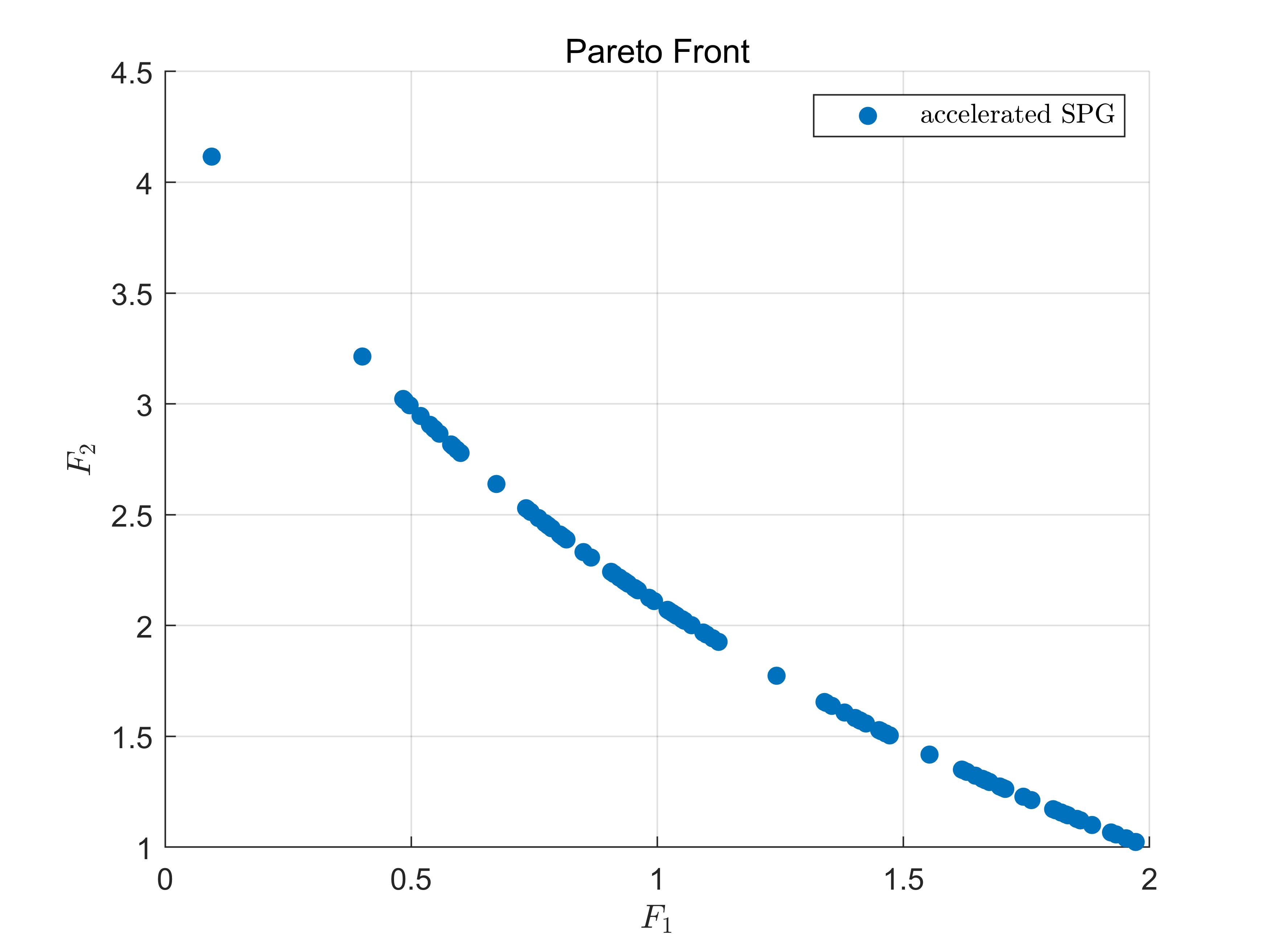}
	\end{minipage}
\begin{minipage}[t] {0.3\textwidth}
	\includegraphics[scale=0.03]{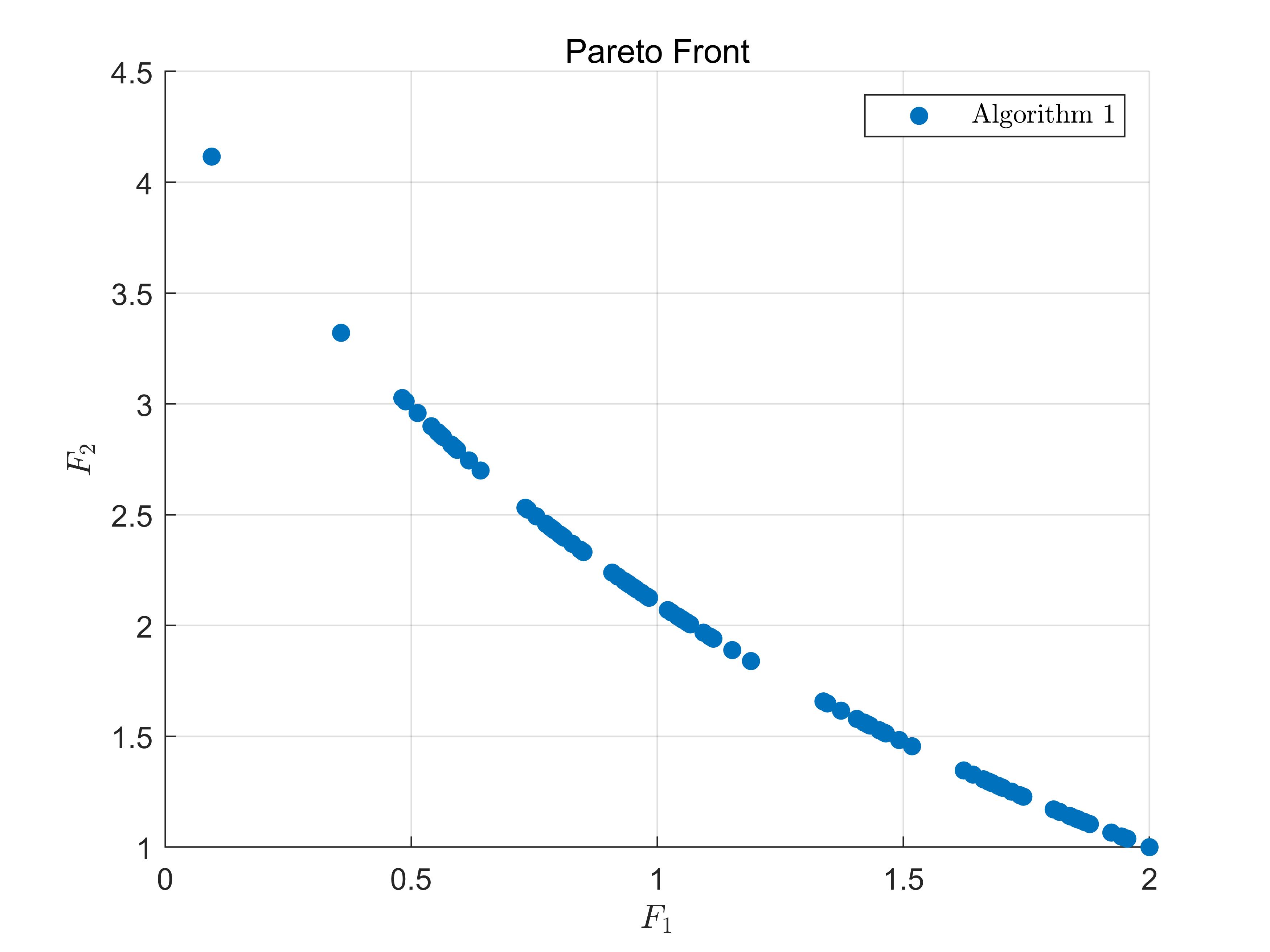}	
\end{minipage}
	\centering\caption{The final Pareto front for JOS1 with $n=500$.}
	\label{fig2}
\end{figure}

\begin{figure}[htpb]
	\centering
	\begin{minipage}[t] {0.4\textwidth}
		\includegraphics[scale=0.04]{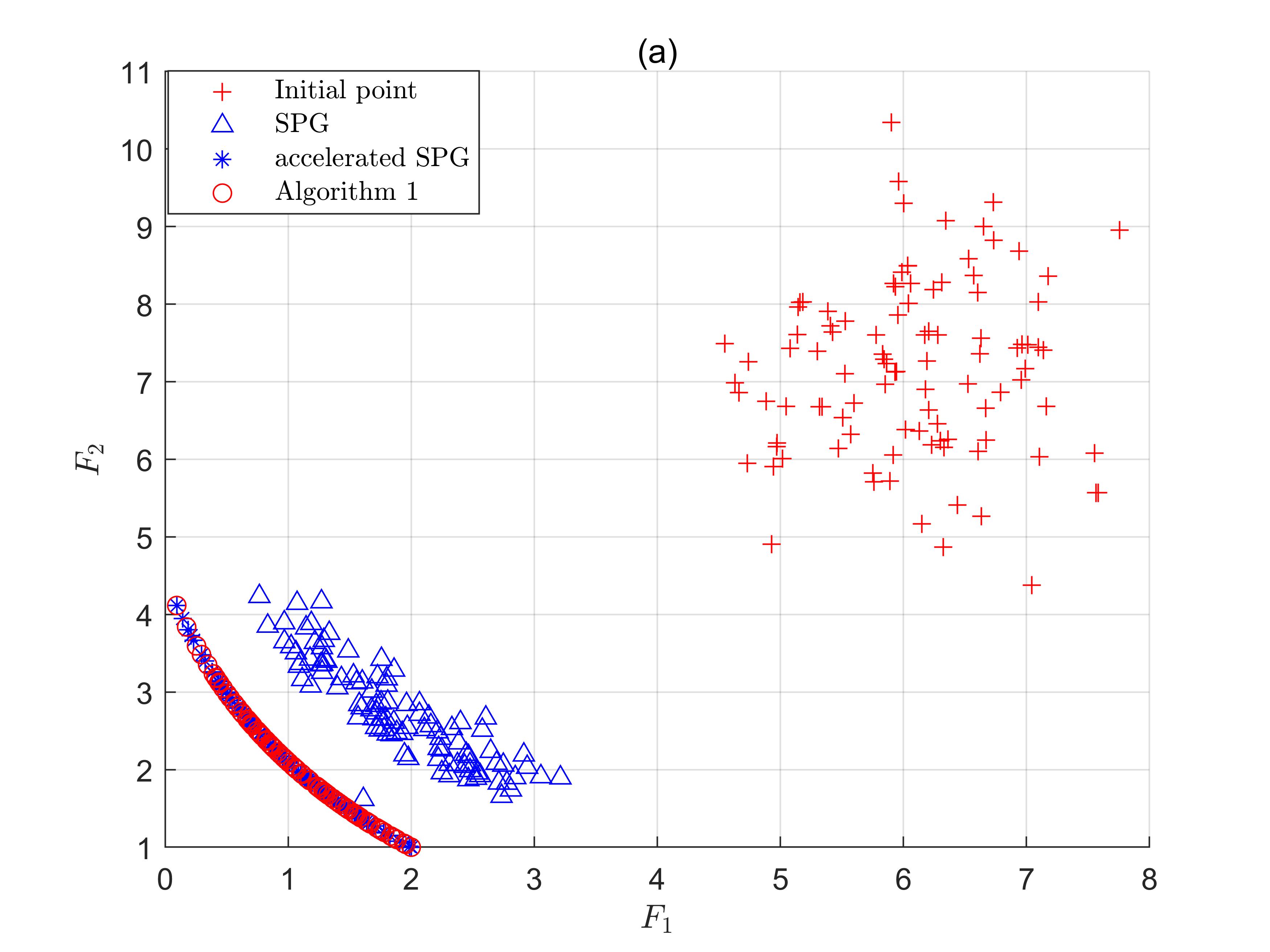}
	\end{minipage}	
	\begin{minipage}[t] {0.4\textwidth}
		\includegraphics[scale=0.04]{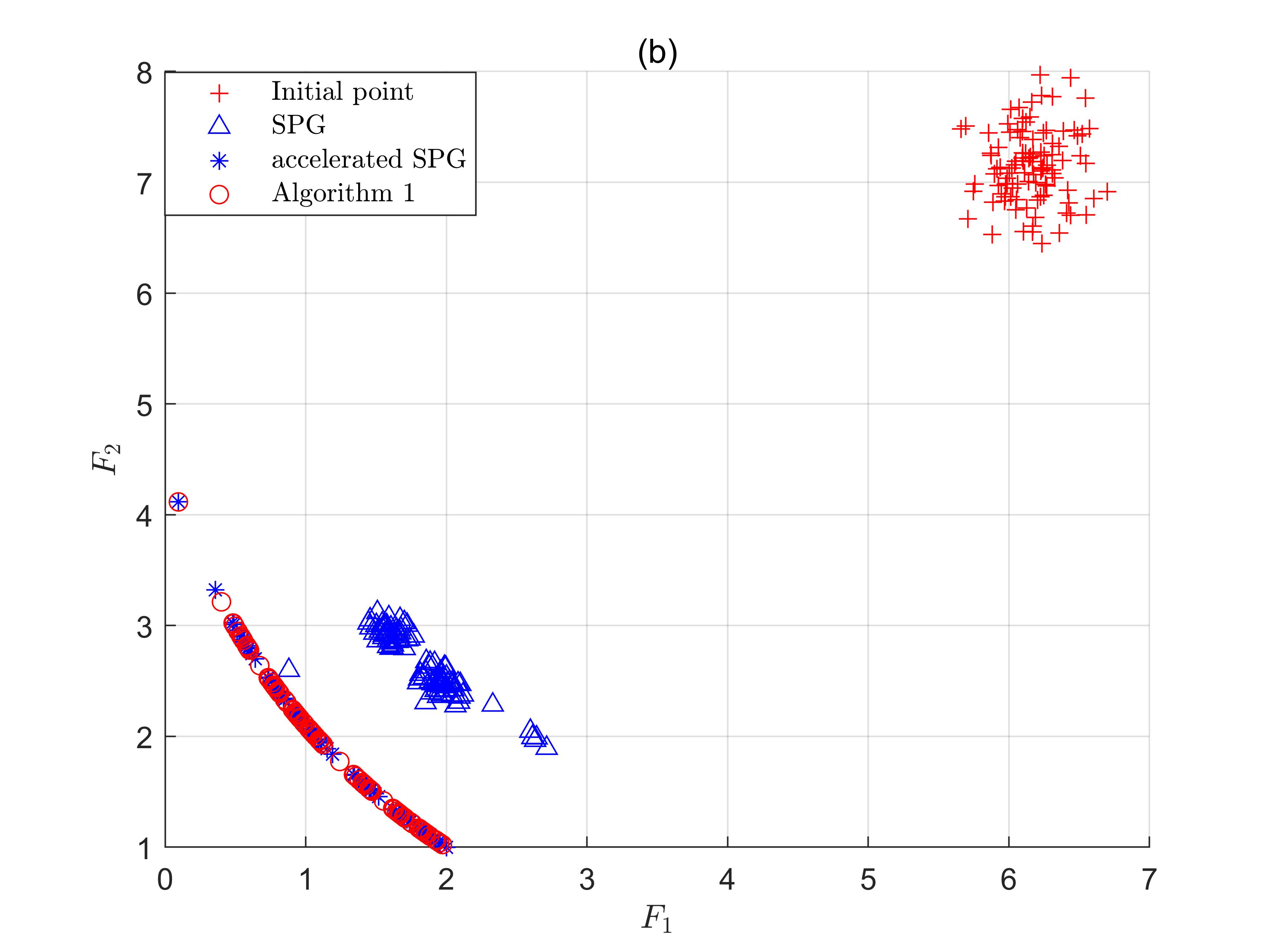}
	\end{minipage}	
	\caption{The  comparison results of Pareto front for JOS1 with
		$n=50$ (left) and $n=500$ (right) when $k=2$.} 
		\label{fig3}
\end{figure}

   Moreover, to illustrate the behaviour of the proposed algorithm,
 we plot the pareto front  for the problem JOS1 with $g_i=l_1$ when $n=50$ and $n=500$. Besides, we make a comparison with the plots of the objective function values for the second iteration points ($k=2$) of  each algorithm under the same initial points, respectively. From Fig. \ref{fig1} to Fig. \ref{fig3}, it can be noted that the pareto front  of the Algorithm \ref{alg1} generated is better than the SPG algorithm.  Besides, the Algorithm \ref{alg1} uses less time than the Algorithm \ref{alg1} and the accelerated SPG.

\begin{table}[htbp]
	\centering
	\caption{Test Problems}
	\begin{tabular}{lclc}
		\toprule
		 {Problem name} & ${m}$ & $n$    & $g_i$ \\
	 	\midrule
		   {JOS1} & 2     & 5,50,500 and 1000 & 0 or $l_1$ \\
		{SD} & 2     & {4} & indicator function \\
		{TOI4} & 2     & {4} &  0 or $l_1$ \\
		{TRIDIA} & 3     & {3} & 0 or $l_1$ \\
		{FDS} & 3     & 5 and 100 & 0 or $l_1$ \\
		\bottomrule
	\end{tabular}%
 \label{table1}
\end{table}%

%

 	\begin{sidewaystable}[htpb]
  		 	\begin{center}
 			\caption{ }\label{tab2}
 			\scalebox{1.0}{
 				\begin{tabular*}{\textheight}{@{\extracolsep{\fill}}lcccccccc@{\extracolsep{\fill}}}
 					\toprule%
 					\multirow{2}{*}{Problem}&\multirow{2}{*}{$n$}&\multirow{2}{*}{$g_i$}& \multicolumn{3}{@{}c@{}}{Averge time($\times 0.01$)}& \multicolumn{3}{@{}c@{}}{Averge iteration counts} 
 					\\\cmidrule{4-6}\cmidrule{7-9}%
 					& & &  {SPG}  & {accelerated SPG}  & Algorithm 1& {SPG} & accelerated SPG &  Algorithm 1\\
 					\midrule
 					
 		{JOS1} & 5  & {0} & 0.080 & 0.070 & 0.070 & 2     & 2     & 2 \\
 		 	 	{JOS1} & 5     & $l_1$    & 0.800 & 0.130 & 0.120 & 20.91 & 2     & 2 \\
 		 	 		{JOS1} & 50  & {0} & 0.120 & 0.110 & 0.090 & 3     & 2     & 2 \\
 		 	 		{JOS1} & 50    & $l_1$    & 1.050 & 0.150 & 0.140 & 22.87 & 2     & 2 \\
 		 	 		{JOS1} & 500   &  {0} & 0.150 & 0.140 & 0.130 & 3     & 2     & 2 \\
 		 	 		{JOS1} & 500   & $l_1$    & 1.200 & 0.250 & 0.220 & 27.08 & 2     & 2 \\
 		 	 		{JOS1} & 1000  & {0} & 0.230 & 0.190 & 0.190 & 3     & 2     & 2 \\
  	 		{JOS1} & 1000  & $l_1$    & 1.330 & 0.280 & 0.280 & 25.48 & 2     & 2 \\
 	        {SD} & 4     & ind   & 10.610 & 6.950 & 5.870 & 885.1 & 827.98 & 777.89 \\
 	 	 {TOI} & 4     & {0} & 0.860 & 0.100 & 0.090 & 1382.14 & 35.47 & 32.36 \\
 	 		 {TOI} & 4     &  $l_1$ & 4.940 & 0.430 & 0.380 & 46.93 & 21.45 & 22.06 \\
 	 		 {TRIDIA} & 3     &  {0} & 1.520 & 1.440 & 1.250 & 1781.86 & 714.53 & 713.46 \\
 	 	 {TRIDIA} & 3     &  $l_1$ & 2.800 & 1.850 & 1.730 & 171.55 & 115.66 & 90.2 \\
 	 	 {FDS} & 5     &  {0} & 1.120 & 0.720 & 0.630 & 151.25 & 148.74 & 132.05 \\
 	 	{FDS} & 5     & $l_1$    & 166.560 & 129.250 & 121.580 & 668.04 & 1037.4 & 1005.38 \\
 	 	 {FDS} & 50    & {0} & 5.010 & 3.070 & 2.800 & 370.24 & 343.09 & 316.12 \\
 	 	{FDS} & 50    & $l_1$   & 490.460 & 390.130 & 365.410 & 1361.48 & 2001  & 2001.000 \\
 	 		 {FDS} & 100   &  {0} & 6.560 & 4.370 & 3.900 & 432.18 & 379.16 & 348.97 \\
 	 	{FDS} & 100   & $l_1$    & 539.420 & 407.550 & 356.840 & 1923.67 & 2001  & 2001 \\

 					\bottomrule
 				\end{tabular*}}
 			  	\end{center}
 		\end{sidewaystable}

	\section{Conclusions} 
	In this paper, we introduce a fast proximal gradient algorithm for multiobjective optimization (denoted as Algorithm 1). We demonstrate that the convergence rate of the accelerated algorithm for multiobjective optimization developed by Tanabe et al. can be enhanced from $O(1/k^2)$ to $o(1/k^2)$ by incorporating a distinct extrapolation term $\frac{k-1}{k+\alpha-1}$ with $\alpha > 3$. 			
	Furthermore, we establish an inexact version of the  Algorithm 1 when the condition of summable error terms is met, and the convergence rate matches that of the Algorithm 1. This resolves the open issues raised in this article \cite{Sonntag1}.
		Finally, we validate the efficiency of the proposed algorithm through some examples.

	\bibliographystyle {plain}
	
	%

	%
	%


\begin{thebibliography}{99}	
		\bibitem{Ansary}
			Ansary, M. A. T. A newton-type proximal gradient method for nonlinear multi-objective
		optimization problems. Optimization Methods and Software, 38:570-590, 2023.
		\bibitem{Attouch}	 
		Attouch, H., Peypouquet, J.: The rate of convergence of Nesterovs accelerated forward-backward 	method is actually faster than $o(1/k^2)$. SIAM J. Optim. 26(3), 1824-1834 (2016)	
			\bibitem{Attouch2}	 		
		Attouch, H., Chbani, Z., Peypouquet, J., Redont, P.: Fast convergence of inertial dynamics and algorithms with asymptotic vanishing viscosity. Math. Program. 168(1-2), 123-175 (2018)	
					\bibitem{Beck}
		 Beck, A.: First-Order Methods in Optimization, Society for Industrial and
		Applied Mathematics, 2017.
		
			
		\bibitem{Carrizo}
		Carrizo, G. A., Lotito, P. A.,  Maciel, M. C. (2016). Trust region globalization strategy for the nonconvex
		unconstrained multiobjective optimization problem. Mathematical Programming, 159, 339-369. 
		
			\bibitem{Chambolle}
		Chambolle, A., Dossal, C.: On the convergence of the iterates of the   ''fast iterative shrinkage/thresholding algorithm'' . J. Optim. Theory Appl. 166, 968-982 (2015)
		
	
		
			\bibitem{ChenJian}
	Chen J, Tang L, Yang X.  A Barzilai-Borwein descent method for multiobjective
	optimization problems. European Journal of Operational Research, 311(1):196-209, 2023.	
		\bibitem{chenw}
	Chen W, Yang X, Zhao Y. Conditional gradient method for vector optimization[J]. Computational Optimization and Applications, 2023: 1-40.
	
			\bibitem{Boyd}	
		Boyd, S., Vandenberghe, L.: Introduction to applied linear algebra: vectors, matrices,
		and least squares. Cambridge university press (2018).
		
		\bibitem{Fliege2}
		Fliege, J.,  Svaiter, B. F. (2000). Steepest descent methods for multicriteria optimization. Mathematical
		Methods of Operations Research, 51, 479-494.
		
	   \bibitem{Fliege1}
		Fliege, J., Drummond, L. M.,  Svaiter, B. F. (2009). Newton's method for multiobjective optimization.
		SIAM Journal on Optimization,20, 602-626.
		
			\bibitem{Fukuda}
		Fukuda E H, Drummond L M G. A survey on multiobjective descent methods[J]. Pesquisa Operacional, 2014, 34: 585-620.
		
		  \bibitem{Lucambio}
		Lucambio Perez, L. R.,  Prudente, L. F. (2018). Nonlinear conjugate gradient methods for vector optimization. SIAM Journal on Optimization, 28, 2690-2720.
		
		
		\bibitem{Tanabe}
		 Tanabe, H., Fukuda, E. H. and Yamashita, N.: New merit functions and error bounds for non-convex multiobjective optimization, arXiv: 2010.09333,
		2020.	
			\bibitem{Tanabe1}
		 Tanabe H, Fukuda E H, Yamashita N. An accelerated proximal gradient method for multiobjective optimization[J]. Computational Optimization and Applications, 2023: 1-35.
		 
		 	\bibitem{Tanabe3}
		 H. Tanabe, E. H. Fukuda, and N. Yamashita. Proximal gradient methods for multiobjective optimization and their applications. Computational Optimization and Applications, 72(2):339-361,
		   2019.
		   	\bibitem{Tanabe4}
		   Tanabe H, Fukuda E H, Yamashita N. A globally convergent fast iterative shrinkage-thresholding algorithm with a new momentum factor for single and multi-objective convex optimization[J]. arXiv preprint arXiv:2205.05262, 2022.
		      	\bibitem{Mita}
		   Mita, K., Fukuda, E.H., Yamashita, N.: Nonmonotone line searches for unconstrained multiobjective
		   optimization problems. J. Global Optim. 75,63-90 (2019)
		   	\bibitem{Nishimura}
		   Nishimura, Yuki, Ellen H. Fukuda, and Nobuo Yamashita. "Monotonicity for Multiobjective Accelerated Proximal Gradient Methods." arXiv preprint arXiv:2206.04412 (2022).
			
	\bibitem{Opial}	
	Opial, Z. (1967). Weak convergence of the sequence of successive approximations for nonexpansive mappings.
		\bibitem{Sonntag1}
	Sonntag K, Peitz S. Fast Multiobjective Gradient Methods with Nesterov Acceleration via Inertial Gradient-like Systems[J]. arXiv preprint arXiv:2207.12707, 2022.
	
	\bibitem{Sonntag2}
	Sonntag K, Peitz S. Fast convergence of inertial multiobjective gradient-like systems with asymptotic vanishing damping[J]. arXiv preprint arXiv:2307.00975, 2023.
	
		
		
	\end{thebibliography}
	
	%
	
\end{document}